\documentclass[english,graybox]{svmult}
\usepackage[T1]{fontenc}
\usepackage[latin9]{inputenc}
\usepackage{babel}
\usepackage{prettyref}
\usepackage{url}
\usepackage{amsmath}
\usepackage{amssymb}
\usepackage{esint}
\usepackage[unicode=true,
 bookmarks=true,bookmarksnumbered=false,bookmarksopen=false,
 breaklinks=false,pdfborder={0 0 1},backref=false,colorlinks=false]
 {hyperref}
\usepackage{breakurl}

\makeatletter
\newenvironment{svmultproof2}{\begin{proof}}{\smartqed\qed\end{proof}}

\usepackage{dimadefs}
\usepackage{mathptmx} 

\DeclareMathOperator{\Op}{\mathfrak{D}}
\DeclareMathOperator{\shift}{E}

\DeclareMathOperator{\id}{I}

\newrefformat{prob}{Problem \ref{#1}}
\newrefformat{prop}{Proposition \ref{#1}}
\newrefformat{def}{Definition \ref{#1}}
\newrefformat{sub}{Section \ref{#1}}
\newrefformat{alg}{Algorithm \ref{#1}}
\newrefformat{con}{Conjecture \ref{#1}}
\newrefformat{cor}{Corollary \ref{#1}}
\newrefformat{rem}{Remark \ref{#1}}
\newrefformat{app}{Appendix \ref{#1}}
\newrefformat{ex}{Example \ref{#1}}
\newrefformat{enu}{Item \ref{#1}}

\makeatother

\begin{document}

\title*{Moment vanishing of piecewise solutions of linear ODEs}

\author{Dmitry Batenkov and Gal Binyamini}

\institute{Dmitry Batenkov \at Department of Mathematics, Weizmann Institute
of Science, Rehovot 76100, Israel. This author is supported by the
Adams Fellowship Program of the Israel Academy of Sciences and Humanities.
\email{dima.batenkov@weizmann.ac.il}\and Gal Binyamini \at Department
of Mathematics, University of Toronto, Canada. \email{galbin@math.utoronto.ca}}

\maketitle

\abstract{We consider the ``moment vanishing problem'' for a general class
of piecewise-analytic functions which satisfy on each continuity interval
a linear ODE with polynomial coefficients. This problem, which essentially
asks how many zero first moments can such a (nonzero) function have,
turns out to be related to several difficult questions in analytic
theory of ODEs (Poincare's Center-Focus problem) as well as in Approximation
Theory and Signal Processing (``Algebraic Sampling''). While the
solution space of any particular ODE admits such a bound, it will
in the most general situation depend on the coefficients of this ODE.
We believe that a good understanding of this dependence may provide
a clue for attacking the problems mentioned above.\\
In this paper we undertake an approach to the moment vanishing problem
which utilizes the fact that the moment sequences under consideration
satisfy a recurrence relation of fixed length, whose coefficients
are polynomials in the index. For any given operator, we prove a general
bound for its moment vanishing index. We also provide uniform bounds
for several operator families.}

\section{Introduction}

\global\long\def\mvi{\ensuremath{\sigma}}
\global\long\def\fun{\ensuremath{f}}
\global\long\def\ord{\ensuremath{n}}
\global\long\def\dg{\ensuremath{d}}
\global\long\def\np{\ensuremath{p}}
\global\long\def\jp{\xi}
\global\long\def\PD{{\cal PD}}
\global\long\def\ff#1#2{\left(#1\right)_{#2}}
\global\long\def\e{\varepsilon}
\global\long\def\mui{\ensuremath{\tau}}
\global\long\def\prp{\mathfrak{R}}

Let $f:\left[a,b\right]\to\reals$ be a bounded piecewise-continuous
function with points of discontinuity (of the first kind)
\[
a=\jp_{0}<\jp_{1}<\dots<\jp_{\np}<\jp_{\np+1}=b,
\]
satisfying on each continuity interval $\left[\jp_{j},\jp_{j+1}\right]$
a linear homogeneous ODE
\begin{equation}
\Op\fun\equiv0,\label{eq:f-satisfy-eqn}
\end{equation}
where $\Op$ is a linear differential operator of order $\ord$ with
polynomial coefficients:
\begin{equation}
\Op=p_{n}\left(x\right)\partial^{n}+\dots+p_{1}\left(x\right)\partial+p_{0}\left(x\right)\id,\qquad\partial=\frac{\dd}{\dd x},\;\deg\; p_{j}\leq\dg_{j}.\label{eq:diffop}
\end{equation}

We say that such $\fun$ belongs to the class $\PD\left(\Op,\np\right)$.
The union of all such $\PD\left(\Op,\np\right)$ is the class $\PD$
of \emph{piecewise D-finite functions}, which was first studied in
\cite{bat2008}.

Any $\fun\in\PD$ has finite moments of all orders:
\begin{equation}
m_{k}\left(\fun\right)=\int_{a}^{b}x^{k}\fun\left(x\right)\dd x,\qquad k=0,1,2,\dots\label{eq:moment-def}
\end{equation}

We consider the following questions.
\begin{problem}
Given $\Op$ and $\np$, determine the \emph{moment vanishing index
}of $\PD\left(\Op,\np\right)$, defined as
\[
\mvi\left(\Op,\np\right)\isdef\sup_{\fun\in\PD\left(\Op,\np\right),\;\fun\not\equiv0}\left\{ k:\; m_{0}\left(\fun\right)=\dots=m_{k}\left(\fun\right)=0\right\} +1.
\]

\end{problem}
In \prettyref{thm:mvi-bound-for-all} below we shall prove that the
moment vanishing index is always finite. Consequently, the following
problem becomes meaningful.
\begin{problem}
\label{prob:uniform-mvi}Find natural families ${\cal F}\subset\PD$
which admit a uniform bound on the moment vanishing index, i.e. for
which
\[
\mvi\left({\cal F}\right)=\sup_{\PD\left(\Op,\np\right)\subset{\cal F}}\mvi\left(\Op,\np\right)<+\infty.
\]

\end{problem}
Our main results, presented in \prettyref{sec:main-results}, provide
a general bound for $\mvi\left(\Op,\np\right)$ in terms of $\Op$.
As a result, several examples of families ${\cal F}$ admitting uniform
bound as above are given. The main technical tool is the recurrence
relation satisfied by the moment sequence, established previously
in \cite{bat2008}.

Our main application is the problem of reconstructing functions $\fun\in\PD$
from a finite number of their moments. Inverse moment problems appear
in some areas of mathematical physics, for instance heat conduction
and inverse potential theory \cite{ang2002mta,gustafsson2000rpd},
as well as in statistics. One particular reconstruction technique,
introduced in \cite{bat2008} and further extended to two-dimensional
setting in \cite{BatGolYom2011}, can be regarded as a prototype for
numerous ``algebraic'' reconstruction methods in signal processing,
such as finite rate of innovation \cite{vetterli2002ssf} and piecewise
Fourier inversion \cite{batFullFourier,batyomAlgFourier}. These methods,
being essentially nonlinear, promise to achieve better reconstruction
accuracy in some cases (as demonstrated recently in \cite{batFullFourier,batyomAlgFourier}),
and therefore we believe their study to be important. In \prettyref{sec:moment-reconstruction}
below we show that an answer to \prettyref{prob:uniform-mvi} would
in turn provide a bound on the minimal number of moments (measurements)
required for unique reconstruction of any $\fun\in{\cal F}$. In essence,
the results of this paper can be regarded as a step towards understanding
the range of applicability of the piecewise D-finite reconstruction
method to general signals in ${\cal PD}$. See \prettyref{sec:moment-reconstruction}
for further details.

Given a family ${\cal F}\subset\PD$, consider the corresponding family
of moment generating functions $\left\{ I_{f}\left(z\right)\right\} _{f\in{\cal F}}$,
where $I_{f}\left(z\right)=\sum_{k=0}^{\infty}m_{k}\left(\fun\right)z^{-k-1}$.
Obtaining information on the moment vanishing index is in fact an
essential step towards studying the analytic properties of $I_{f}$,
in particular a bound on its number of zeros near infinity (as provided
by the notion of ``Taylor Domination'', see \cite{batyomTaylrPoincBarc}),
as well as conditions for its identical vanishing. In turn, these
questions play a central role in studies of the Center-Focus and Smale-Pugh
problems for the Abel differential equation, see \cite{briskin2010center,pakovich:cti}
and references therein.

The moment vanishing problem has been previously studied in the complex
setting by V.Kisunko \cite{kisunko1}. He showed that a uniform bound
$\mvi\left({\cal F}\right)$ exists for families ${\cal F}$ consisting
of non-singular operators, by using properties of Cauchy type integrals.
In contrast, in this paper we consider the real setting only, while
proving uniform bounds for some singular (as well as regular) operator
families. Our method is based on the linear recurrence relation satisfied
by the moment sequence. Using this method, in \prettyref{sec:mgf}
we provide an alternative proof of Kisunko's result, stating that
the moment generating function $I_{\fun}\left(z\right)$ of some $\fun\in\PD\left(\Op,\np\right)$
satisfies a non-homogeneous ODE
\[
\Op I_{\fun}\left(z\right)=R_{f}\left(z\right)
\]
for a very special rational function $R_{f}\left(z\right)$, which
depends on $\Op$ and on the values of $\fun$ at the discontinuities.

In \prettyref{sec:fuchsian} we provide an interpretation of our main
result in the language of Fuchsian theory of ODE.

\subsection{Acknowledgements}

The authors would like to thank Y.Yomdin for useful discussions.

\section{\label{sec:moment-reconstruction}Moment reconstruction}

We start by defining some preliminary notions.
\begin{definition}
The Pochhammer symbol $\ff ij$ denotes the falling factorial 
\[
\ff ij\isdef i(i-1)\cdot\dotsc\cdot(i-j+1),\qquad i\in\reals,\; j\in\naturals
\]
and the expression $\ff ij$ is defined to be zero for $i<j$.
\end{definition}
\begin{minipage}[t]{1\columnwidth}%
\end{minipage}
\begin{definition}
Given $\Op$ of the form \eqref{eq:diffop}, the \emph{bilinear concomitant}
(\cite[p.211]{ince1956ode}) is the homogeneous bilinear form, defined
for any pair of sufficiently smooth functions $u\left(x\right),v\left(x\right)$
as follows (all symbols depend on $x$):
\begin{eqnarray}
P_{\Op}\left(u,v\right) & \isdef & u\left\{ p_{1}v-\partial\left(p_{2}v\right)+\dots+\left(-1\right)^{\ord-1}\partial^{\ord-1}\left(p_{\ord}v\right)\right\} \label{eq:bilinear-concomitant}\\
 & + & u'\left\{ p_{2}v-\partial\left(p_{3}v\right)+\dots+\left(-1\right)^{\ord-2}\partial^{\ord-2}\left(p_{\ord}v\right)\right\} \nonumber \\
 & + & \dots\nonumber \\
 & + & u^{\left(n-1\right)}\cdot\left(p_{\ord}v\right).\nonumber 
\end{eqnarray}
\end{definition}
\begin{proposition}[Green's formula, \cite{ince1956ode}]
Given $\Op$ of the form \eqref{eq:diffop}, let the formal adjoint
operator be defined by
\[
\Op^{*}\left\{ \cdot\right\} \isdef\sum_{j=0}^{n}\left(-1\right)^{j}\partial^{j}\left\{ p_{j}\left(x\right)\cdot\right\} .
\]
Then for any pair of sufficiently smooth functions $u\left(x\right),v\left(x\right)$
the following identity holds:
\begin{equation}
\int_{a}^{b}v\left(x\right)\left(\Op u\right)\left(x\right)\D x-\int_{a}^{b}u\left(x\right)\left(\Op^{*}v\right)\left(x\right)\D x=P_{\Op}\left(u,v\right)\left(b\right)-P_{\Op}\left(u,v\right)\left(a\right).\label{eq:green-formula}
\end{equation}
\end{proposition}
\begin{theorem}[\cite{bat2008}]
\label{thm:dfinite-inv}Let $\fun\in{\cal PD}\left(\Op,\np\right)$
with $\Op$ of the form \eqref{eq:diffop}. Denote the discontinuities
of $\fun$ by $a=\jp_{0}<\jp_{1}<\dots<\jp_{\np}<\jp_{\np+1}=b$.
Then the moments $m_{k}=\int_{a}^{b}\fun\left(x\right)\D x$ satisfy%
\footnote{For consistency of notation, the sequence $\left\{ m_{k}\right\} $
is understood to be extended with zeros for negative $k$.%
} the recurrence relation
\begin{equation}
\sum_{j=0}^{n}\sum_{i=0}^{d_{j}}a_{i,j}\left(-1\right)^{j}\ff{i+k}jm_{i-j+k}=\e_{k},\quad k=0,1,\dots,\label{eq:moments-main-rec}
\end{equation}
where 
\begin{equation}
\e_{k}=-\sum_{j=0}^{\np}\left\{ P_{\Op}\left(\fun,x^{k}\right)\left(\jp_{j+1}^{-}\right)-P_{\Op}\left(\fun,x^{k}\right)\left(\jp_{j}^{+}\right)\right\} .\label{eq:ek-def}
\end{equation}
\end{theorem}
\begin{svmultproof2}
Apply Green's formula \eqref{eq:green-formula} to the identity
\[
\int_{\jp_{j}}^{\jp_{j+1}}x^{k}\left(\Op f\right)\left(x\right)\D x\equiv0
\]
for each $j=0,\dots,\np$ and sum up. The result is
\begin{eqnarray*}
\sum_{j=0}^{\np}\int_{\jp_{j}}^{\jp_{j}+1}\fun\left(x\right)\Op^{*}\left\{ x^{k}\right\} \D x & = & -\sum_{j=0}^{\np}\left\{ P_{\Op}\left(\fun,x^{k}\right)\left(\jp_{j+1}^{-}\right)-P_{\Op}\left(\fun,x^{k}\right)\left(\jp_{j}^{+}\right)\right\} \\
\int_{a}^{b}\fun\left(x\right)\Op^{*}\left\{ x^{k}\right\} \D x & = & \e_{k}
\end{eqnarray*}
The left-hand side of the last formula is precisely the linear combination
of the moments given by the left-hand side of \eqref{eq:moments-main-rec}.
This finishes the proof.
\end{svmultproof2}
Now consider the problem of recovering $\fun\in\PD\left(\Op,\np\right)\subset\PD$
from the moments $\left\{ m_{0}\left(\fun\right),\dots,m_{N}\left(\fun\right)\right\} $
(the operator $\Op$ is assumed unknown in the most general setting).
Based on the recurrence relation \eqref{eq:moments-main-rec}, we
demonstrate in \cite{bat2008} that an exact recovery is possible,
provided that the number $N\in\naturals$ is sufficiently large. However,
the question of obtaining an upper bound for $N$ turns out to be
non-trivial, as we now demonstrate.
\begin{definition}
Given $\Op$ and $\np$, the \emph{moment uniqueness index} $\mui\left(\Op,\np\right)$
is defined by
\[
\mui\left(\Op,\np\right)\isdef\sup_{f,g\in\PD\left(\Op,\np\right),\; f\not\equiv g}\left\{ k:\; m_{j}\left(f\right)=m_{j}\left(g\right),\;0\leqslant j\leqslant k\right\} +1.
\]

\end{definition}
In other words, given $\Op$ and $\np$, at least $\mui\left(\Op,\np\right)$
first moments of $\fun\in\PD\left(\Op,\np\right)$ are necessary for
unique reconstruction of $\fun$.

Recalling boundedness of $\mvi\left(\Op,\np\right)$ (see \prettyref{thm:mvi-bound-for-all}
below), we immediately obtain the following conclusion.
\begin{lemma}
For any operator $\Op$ and any $\np$ 
\[
\mui\left(\Op,\np\right)\leq\mvi\left(\Op,2\np\right).
\]
\end{lemma}
\begin{svmultproof2}
Let $N=\mvi\left(\Op,2\np\right)$. Take$\fun_{1},\fun_{2}$ having
$\np$ jump points each, satisfying $\Op\fun_{1}\equiv0,\Op\fun_{2}\equiv0$
on each continuity interval such that
\begin{eqnarray*}
m_{0}\left(\fun_{1}\right) & = & m_{0}\left(\fun_{2}\right)\\
 & \dots\\
m_{N}\left(\fun_{_{1}}\right) & = & m_{N}\left(\fun_{2}\right).
\end{eqnarray*}
The function $g=\fun_{1}-\fun_{2}$ has at most $2\np$ jumps, and
it satisfies $\Op g\equiv0$ on each continuity interval. The first
$N$ moments of $g$ are zero, therefore $g\equiv0$ and thus $\fun_{1}\equiv\fun_{2}$.
Therefore $\mui\left(\Op,\np\right)\leq N$.
\end{svmultproof2}
Consequently, in order to uniquely reconstruct an unknown $\fun\in{\cal F}\subset\PD$,
it is sufficient to get a uniform bound $\mvi\left({\cal F}\right)$$ $
for the family ${\cal F}$. Perhaps the most natural choice for such
families is when the parameters $\np,n,\left\{ d_{j}\right\} _{j=0}^{n}$
are fixed. Unfortunately, without making additional assumptions, the
moment vanishing index of such families cannot be uniformly bounded.
This can be seen from the following example.
\begin{example}
\label{ex:legendre}Let $\Op_{m}$ denote the Legendre differential
operator
\[
\Op_{m}=\left(1-x^{2}\right)\frac{\dd^{2}}{\dd x^{2}}-2x\frac{\dd}{\dd x}+m\left(m+1\right)\id,
\]
and consider $a=-1,\; b=1$ and $\np=0$. It is well-known that for
each $m\in\naturals$, the regular solution of $\Op_{m}\fun=0$ is
${\cal L}_{m}$ - the Legendre polynomial of degree $m$. Since the
first $m-1$ moments of ${\cal L}_{m}$ are zero, we conclude that
\[
\mvi\left(\Op_{m}\right)=m
\]
and therefore $\mvi\left(\Op\right)$ cannot be uniformly bounded
in terms of the combinatorial type of $\Op$ only.
\end{example}
Using the subsequent results, in \prettyref{sec:fuchsian} we shall
in fact provide an explanation of this behaviour.

\section{\label{sec:powersums}Generalized power sums}
\begin{proposition}
\label{prop:ek-structure}The sequence $\left\{ \e_{k}\right\} $,
given by \prettyref{thm:dfinite-inv}, is of the form
\begin{equation}
\e_{k}=\sum_{j=0}^{\np+1}\sum_{\ell=0}^{n-1}\jp_{j}^{k-\ell}\ff k{\ell}c_{\ell,j},\label{eq:ek-as-exp-poly}
\end{equation}
where each $c_{\ell,j}$ is a homogeneous bilinear form in the two
sets of variables
\begin{align*}
 & \{p_{m}(\jp_{j}),p'_{m}(\jp_{j}),\dots,p_{m}^{(n-1)}(\jp_{j})\}_{m=0}^{n},\\
 & \{\fun(\jp_{j}^{+})-\fun(\jp_{j}^{-}),\;\fun'(\jp_{j}^{+})-\fun'(\jp_{j}^{-}),\dots,\fun^{(n-1)}(\jp_{j}^{+})-\fun^{(n-1)}(\jp_{j}^{-})\}.
\end{align*}
\end{proposition}
\begin{svmultproof2}
Denote for convenience $\fun\left(a^{-}\right)=f\left(b^{+}\right)=0$.
Now consider the definition of $\left\{ \e_{k}\right\} $ given by
\eqref{eq:ek-def}. Rearranging terms, we write
\[
\e_{k}=\sum_{j=0}^{\np+1}\left\{ P_{\Op}\left(\fun,x^{k}\right)\left(\jp_{j}^{+}\right)-P_{\Op}\left(\fun,x^{k}\right)\left(\jp_{j}^{-}\right)\right\} .
\]
Furthermore, using the fact that the functions $\left\{ p_{m}\left(x\right)\right\} _{m=0}^{n}$
and $x^{k}$ are continuous at each $\jp_{j}$, we have{\small 
\begin{eqnarray}
P_{\Op}\left(f,x^{k}\right)\left(\jp_{j}^{+}\right)-P_{\Op}\left(\fun,x^{k}\right)\left(\jp_{j}^{-}\right) & =\label{eq:ek-term-j}\\
 & = & \left\{ \fun\left(\jp_{j}^{+}\right)-\fun\left(\jp_{j}^{-}\right)\right\} \times\nonumber \\
 &  & \;\times\left\{ p_{1}\left(\jp_{j}\right)\jp_{j}^{k}-\left(p_{2}\left(\jp_{j}\right)k\jp_{j}^{k-1}+p_{2}'\left(\jp_{j}\right)\jp_{j}^{k}\right)+\dots\right\} \nonumber \\
 &  & +\dots\nonumber \\
 &  & +\left\{ \fun^{\left(n-1\right)}\left(\jp_{j}^{+}\right)-\fun^{\left(n-1\right)}\left(\jp_{j}^{-}\right)\right\} p_{n}\left(\jp_{j}\right)\jp_{j}^{k}.\nonumber 
\end{eqnarray}
}Now using the definition \eqref{eq:bilinear-concomitant}, the claim
is evident.
\end{svmultproof2}
The expression \eqref{eq:ek-as-exp-poly} for $\e_{k}$ is nothing
else but a generalized power sum. Let us recall several well-known
facts about them (see e.g. \cite[Section 2.3]{elaydi2005ide} or \cite{myerPort1995}).
\begin{proposition}
Let the sequence $s_{k}$ be of the form 
\begin{equation}
s_{k}=\sum_{j=0}^{\np+1}\sum_{\ell=0}^{\ord-1}a_{\ell,j}\ff k{\ell}\jp_{j}^{k-\ell}\qquad a_{\ell,j},\jp_{j}\in\complexfield.\label{eq:conf-prony}
\end{equation}
Then it satisfies the following linear recurrence relation with constant
coefficients of length $\ord\left(\np+2\right)+1$:
\begin{align}
\biggl(\prod_{j=0}^{\np+1}(\shift-\jp_{j}\id)^{\ord}\biggr)s_{k} & =0\label{eq:const-coef-rec}
\end{align}
 where $\shift$ is the forward shift operator in $k$ and $\id$
is the identity operator.

Conversely, the fundamental set of solutions of the recurrence relation
\eqref{eq:const-coef-rec} is
\[
\left\{ \jp_{0}^{k},k\jp_{0}^{k-1},\dots,\ff k{n-1}\jp_{0}^{k-n+1},\;\dots,\;\jp_{\np+1}^{k},k\jp_{\np+1}^{k-1},\dots,\ff k{n-1}\jp_{\np+1}^{k-n+1}\right\} .
\]
\end{proposition}
\begin{corollary}
\label{cor:max-zeros-conf-prony}The sequence $s_{k}$ as above, which
is not identically zero, can have at most $\ord\left(\np+2\right)-1$
first consecutive zero terms $s_{0}=\dots=s_{n\left(\np+2\right)-2}=0$.\end{corollary}
\begin{svmultproof2}
If $s_{0}=\dots=s_{n\left(\np+2\right)-1}=0,$ then by the recurrence
relation \eqref{eq:const-coef-rec} we would have automatically $s_{n\left(\np+2\right)}=s_{n\left(\np+2\right)+1}=\dots=0$.
\end{svmultproof2}
\begin{minipage}[t]{1\columnwidth}%
\end{minipage}
\begin{corollary}
\label{cor:vanishing-ek-implies-vanishing-aij}Assume that the numbers
$\left\{ \jp_{j}\right\} _{j=0}^{\np+1}\subset\complexfield$ are
pairwise distinct. Let the sequence $s_{k}$ be given by \eqref{eq:conf-prony},
with a-priori unknown $\left\{ a_{i,j}\right\} $. If $s_{k}=0$ for
all $k\in\naturals,$ then necessarily all the coefficients $\left\{ a_{i,j}\right\} $
are zero.
\end{corollary}

\section{\label{sec:main-results}Main results}

Let us now return to our main goal, namely, obtaining upper bounds
on the moment vanishing index $\mvi\left(\Op,\np\right)$ . 
\begin{definition}
\label{def:alphas}Given $\Op$ of the form \eqref{eq:diffop}, denote
for each $j=0,\dots,n$
\[
\alpha_{j}\isdef d_{j}-j,
\]
and also
\[
\alpha=\alpha\left(\Op\right)\isdef\max_{j=0,\dots,n}\alpha_{j}.
\]
\end{definition}
\begin{proposition}
\label{prop:mom-van-implies-ek-van}Let $\fun\in\PD\left(\Op,\np\right)$.
Then vanishing of the first $\left(\np+2\right)\ord+\alpha\left(\Op\right)$
moments of $\fun$ (i.e. $m_{0}=\dots=m_{\left(\np+2\right)n+\alpha-1}=0$)
implies identical vanishing of the sequence $\left\{ \e_{k}\right\} $
defined by \prettyref{thm:dfinite-inv}.\end{proposition}
\begin{svmultproof2}
Consider the recurrence relation \eqref{eq:moments-main-rec}. Denote
its left-hand side by $\mu_{k}$. Obviously, since each $\mu_{k}$
is a linear combination of the moments, we have
\[
\mu_{0}=\dots=\mu_{n\left(\np+2\right)-1}=0.
\]
Consequently, the corresponding right-hand sides also vanish, i.e.
\begin{equation}
\e_{0}=\dots=\e_{n\left(\np+2\right)-1}=0.\label{eq:ek-vanish}
\end{equation}
The conclusion follows immediately from \prettyref{cor:max-zeros-conf-prony}.
\end{svmultproof2}
Now we establish our main result.
\begin{theorem}
\label{thm:mvi-general-case}Let $\fun\in\PD\left(\Op,\np\right)$,
$\fun\not\equiv0$ with discontinuity points
\[
a=\jp_{0}<\jp_{1}<\dots<\jp_{\np}<\jp_{\np+1}=b.
\]
Assume that $p_{n}\left(\jp_{j}\right)\neq0$ for at least one $\jp_{j}$
as above. Then at most
\[
\left(\np+2\right)\ord+\alpha\left(\Op\right)-1
\]
first moments of $\fun$ can vanish (i.e. $m_{0}=\dots=m_{\left(\np+2\right)n+\alpha-2}=0$).\end{theorem}
\begin{svmultproof2}
Assume by contradiction that the first $\left(\np+2\right)\ord+\alpha$
moments of $\fun$ vanish, i.e.
\[
m_{0}=\dots=m_{n\left(\np+2\right)+\alpha-1}=0.
\]
By \prettyref{prop:mom-van-implies-ek-van} and \prettyref{cor:vanishing-ek-implies-vanishing-aij}
we immediately conclude that
\[
c_{\ell,j}=0,\quad j=0,\dots,\np+1,\;\ell=0,\dots,n-1,
\]
where $\left\{ c_{\ell,j}\right\} $ are described by \prettyref{prop:ek-structure}.
Now we take the concrete $j$ for which $p_{n}\left(\jp_{j}\right)\neq0$.
This means that the operator $\Op$ is regular at $\jp_{j}$, and
consequently each solution to $\Op\fun=0$ in the neighborhood of
$\jp_{j}$ is uniquely determined by the initial values $\fun\left(\jp_{j}\right),\dots,\fun^{\left(n-1\right)}\left(\jp_{j}\right)$.
We claim that
\begin{equation}
\fun\left(\jp_{j}^{+}\right)-\fun\left(\jp_{j}^{-}\right)=\fun'\left(\jp_{j}^{+}\right)-\fun'\left(\jp_{j}^{-}\right)=\dots=\fun^{\left(n-1\right)}\left(\jp_{j}^{+}\right)-\fun^{\left(n-1\right)}\left(\jp_{j}^{-}\right)=0.\label{eq:vanishing-jumps}
\end{equation}
In this case, we would immediately conclude that the function $\fun$
is analytic at $\jp_{j}$ (being a solution of analytic ODE), contradicting
the assumption that $\jp_{j}$ is a point of discontinuity of $\fun$.

To prove \eqref{eq:vanishing-jumps}, we proceed as follows. By \prettyref{prop:ek-structure}
it is easy to see that the term $c_{n-1,j}\ff k{n-1}\jp_{j}^{k-n+1}$
is in fact equal to
\[
\left\{ \fun\left(\jp_{j}^{+}\right)-\fun\left(\jp_{j}^{-}\right)\right\} \ff k{n-1}p_{n}\left(\jp_{j}\right)
\]
in the expression for $\e_{k}$. Since $p_{n}\left(\jp_{j}\right)\neq0$,
we conclude that $\fun\left(\jp_{j}^{+}\right)-\fun\left(\jp_{j}^{-}\right)=0$.
Substituting this into \eqref{eq:ek-term-j}, we see that the next
term $c_{n-2,j}\ff k{n-2}\jp_{j}^{k-n+2}$ equals
\[
\left\{ \fun'\left(\jp_{j}^{+}\right)-\fun'\left(\jp_{j}^{-}\right)\right\} \ff k{n-2}\jp_{j}^{k-n+2}p_{n}\left(\jp_{j}\right),
\]
and thus $\fun'\left(\jp_{j}^{+}\right)-\fun'\left(\jp_{j}^{-}\right)=0$.
Proceeding in this manner, we arrive at \eqref{eq:vanishing-jumps}.
This finishes the proof of \prettyref{thm:mvi-general-case}.
\end{svmultproof2}
As a first consequence, we have the real-valued version of the result
by Kisunko \cite{kisunko1}.
\begin{corollary}
For every $n,d>0$ and $\np\geqslant0$ consider the family 
\[
{\cal F}_{n,\np,d}^{\left(1\right)}=\left\{ \fun\in\PD\left(\Op,\np\right):\quad\Op=\sum_{j=0}^{n}p_{j}\left(x\right)\partial^{j},\;\alpha\left(\Op\right)=d,\; p_{n}\left(x\right)\neq0\;\text{on }\left[a,b\right]\right\} .
\]
Then
\[
\mvi\left({\cal F}_{n,\np,d}^{\left(1\right)}\right)\leqslant\left(\np+2\right)n+d-1.
\]

\end{corollary}
Since the leading coefficient $p_{n}\left(x\right)$ cannot vanish
at more than $\deg p_{n}$ points, we also have the following result.
\begin{corollary}
For every $n,d>0$ and $p\geqslant0$ consider the family
\[
{\cal F}_{n,\np,d}^{\left(2\right)}=\left\{ \fun\in\PD\left(\Op,\np\right):\quad\Op=\sum_{j=0}^{n}p_{j}\left(x\right)\partial^{j},\quad\alpha\left(\Op\right)=d,\;\deg p_{n}<\np+2\right\} .
\]
Then
\[
\mvi\left({\cal F}_{n,\np,d}^{\left(2\right)}\right)\leqslant\left(\np+2\right)n+d-1.
\]

\end{corollary}
Let us now try to establish what happens in the general case. Let
$\fun\in\PD\left(\Op,\np\right)$, $\fun\not\equiv0$. Consider two
possibilities.
\begin{enumerate}
\item The sequence $\left\{ \e_{k}\right\} $ does not vanish identically.
In this case, at least some of its initial terms $\left\{ \e_{0},\dots,\e_{n\left(\np+2\right)-1}\right\} $
must be nonzero (\prettyref{cor:max-zeros-conf-prony}). But this
immediately implies that some of the first $n\left(\np+2\right)+\alpha-1$
moments must be nonzero as well (otherwise the equality \eqref{eq:moments-main-rec}
cannot hold).
\item The sequence $\left\{ \e_{k}\right\} $ vanishes identically, but
\prettyref{thm:mvi-general-case} is not applicable. In this case
the recurrence relation \eqref{eq:moments-main-rec} becomes homogeneous.
We rewrite it in the form
\begin{equation}
\sum_{\ell=-n}^{\alpha}q_{\ell}\left(k\right)m_{k+\ell}=0,\qquad k=0,1,\dots,\label{eq:muk-homo}
\end{equation}
where
\begin{equation}
q_{\ell}\left(k\right)\isdef\sum_{j=0}^{n}\left(-1\right)^{j}a_{\ell+j,j}\ff{k+\ell+j}j.\label{eq:ql-def}
\end{equation}
The leading coefficient $q_{\alpha}\left(k\right)$ may have positive
integer zeros. Let $\Lambda\left(\Op\right)$ denote the largest such
zero. Then we claim that no more than $\alpha+\Lambda\left(\Op\right)$
moments can vanish. Indeed, starting with $k=\Lambda\left(\Op\right)+1$
we can safely divide the recurrence \eqref{eq:muk-homo} by $q_{\alpha}\left(k\right)$
and obtain
\[
m_{k+\alpha}=\sum_{\ell=-n}^{\alpha-1}r_{\ell}\left(k\right)m_{k+\ell},\quad k\geqslant\Lambda\left(\Op\right)+1,
\]
where $r_{\ell}\left(k\right)$ are some rational functions with non-vanishing
denominators. Therefore if the first $\Lambda\left(\Op\right)+\alpha+1$
moments are zero, then all the rest of the moments must vanish, implying
vanishing of $\fun$ itself.
\end{enumerate}
Thus we have proved the following result.
\begin{theorem}
\label{thm:mvi-bound-for-all}For every $\Op,\np$ we have
\[
\mvi\left(\Op,\np\right)\leqslant\max\left\{ n\left(\np+2\right)-1,\;\Lambda\left(\Op\right)\right\} +\alpha\left(\Op\right).
\]

\end{theorem}
In \prettyref{sec:fuchsian}, we demonstrate that in the case of Fuchsian
differential operators, the number $\Lambda\left(\Op\right)$ has
a well-known interpretation.

\section{\label{sec:mgf}Moment generating function}

In this section we provide an alternative proof for the result of
Kisunko \cite{kisunko1} concerning moment generating functions.
\begin{proposition}
\textup{Let $\fun\in\PD$.} The formal power series
\[
I_{\fun}\left(z\right)\isdef\sum_{k=0}^{\infty}\frac{m_{k}}{z^{k+1}}
\]
is in fact the Laurent series of the Cauchy type integral \textup{
\[
\int_{a}^{b}\frac{\fun\left(t\right)\dd t}{z-t}.
\]
}\end{proposition}
\begin{svmultproof2}
Write $\frac{1}{z-t}=\frac{1}{z}\left(\frac{1}{1-\frac{t}{z}}\right)$
and expand into geometric series. Convergence follows immediately
for $z\to\infty$.
\end{svmultproof2}
The generalized power sums (\prettyref{sec:powersums}) also have
a well-known interpretation as the Taylor coefficients of rational
functions. The following fact is well-known, and so we omit the proof.
\begin{proposition}
\label{prop:rational-gen-fun}Let the sequence $\left\{ s_{k}\right\} $
be of the form \eqref{eq:conf-prony}. Then the formal generating
function
\[
g\left(z\right)=\sum_{k=0}^{\infty}\frac{s_{k}}{z^{k+1}}
\]
is a regular at infinity rational function, with poles $\left\{ \jp_{0},\dots,\jp_{\np+1}\right\} $,
each with multiplicity at most $n$. In particular,
\begin{equation}
g\left(z\right)=\sum_{j=0}^{\np+1}\sum_{\ell=0}^{n-1}\frac{\left(-1\right)^{\ell}\ell!a_{\ell,j}}{\left(z-\jp_{j}\right)^{\ell+1}}.\label{eq:rational-fn-explicit}
\end{equation}
\end{proposition}
\begin{theorem}
\label{thm:cauchy-int-ode}Let $\fun\in\PD\left(\Op,\np\right)$.
Then the Cauchy integral $I_{\fun}$ satisfies in the neighborhood
of $\infty$ the inhomogeneous ODE 
\begin{equation}
\Op I_{\fun}\left(z\right)=R_{\fun}\left(z\right),\label{eq:cauchy-int-ode}
\end{equation}
where $R_{\fun}\left(z\right)$ is the rational function whose Taylor
coefficients at infinity are given by the sequence $\epsilon_{k}$
as in \eqref{eq:ek-def}. Consequently, $R_{f}\left(z\right)$ is
given by the explicit expression \eqref{eq:rational-fn-explicit},
with $a_{\ell,j}$ replaced by $c_{\ell,j}$ from \eqref{eq:ek-as-exp-poly}
(\prettyref{prop:ek-structure}).\end{theorem}
\begin{svmultproof2}
Consider the asymptotic expansion of the function $\Op I_{\fun}$
at infinity
\[
\Op I_{\fun}=\sum_{k=0}^{\infty}\frac{s_{k}}{z^{k+1}}.
\]
By substituting $\Op$ as in \eqref{eq:diffop} and $I_{f}=\sum_{k=0}^{\infty}\frac{m_{k}}{z^{k+1}}$
we get
\begin{eqnarray*}
\Op I_{\fun} & = & \sum_{j=0}^{n}p_{j}\left(z\right)I_{\fun}^{\left(j\right)}\left(z\right)=\sum_{j=0}^{n}\sum_{i=0}^{d_{j}}\sum_{k=0}^{\infty}m_{k}\frac{\left(-1\right)^{j}\ff{k+j}j}{z^{k+1+j-i}}a_{i,j}\\
_{\left(k+j-i\to t\right)} & = & \sum_{t=0}^{\infty}\frac{1}{z^{t+1}}\sum_{j=0}^{n}\left(-1\right)^{j}a_{i,j}\ff{t+i}jm_{t+i-j}\\
 & = & \sum_{s=0}^{\infty}\frac{s_{k}}{z^{k+1}}.
\end{eqnarray*}
Comparing powers of $z$ we have that $s_{k}=\mu_{k}$ where $\mu_{k}$
denote the left-hand side of \eqref{eq:moments-main-rec}. From $\mu_{k}=\epsilon_{k}$
the conclusion follows.
\end{svmultproof2}

\section{\label{sec:fuchsian}Fuchsian operators}

In this section we employ notions from the classical Fuchsian theory
of linear ODEs in the complex domain (we used the reference \cite{henrici1977applied}).

Assume that the sequence $\left\{ \e_{k}\right\} $ vanishes identically.
In this case, the Cauchy integral $I_{f}$ satisfies in the neighborhood
of $\infty$ the \emph{homogeneous ODE
\[
\Op I_{f}=0.
\]
}
\begin{definition}
The operator $\Op$ is said to belong to the class $\prp$ if it has
at most a regular singularity at $\infty$.\end{definition}
\begin{lemma}
Let $\Op\in\prp$. Then
\begin{enumerate}
\item The numbers $\alpha_{j}$ (see \prettyref{def:alphas}) satisfy
\begin{equation}
\alpha_{n}\geqslant\alpha_{j},\quad j=0,\dots,n-1.\label{eq:fuchs-cond}
\end{equation}

\item The characteristic exponents of $\Op$ at the point $\infty$ are
the roots of the equation
\[
q_{\alpha_{n}}\left(s-1\right)=0,
\]
where $q_{\ell}\left(k\right)$ is defined by \eqref{eq:ql-def}.
\end{enumerate}
\end{lemma}
\begin{svmultproof2}
Dividing the coefficients of $\Op$ by $p_{n}$, we get the operator
\[
\partial^{n}+r_{1}\left(z\right)\partial^{n-1}+\dots+r_{n}\left(z\right)\id,\qquad r_{j}\left(z\right)=\frac{p_{n-j}\left(z\right)}{p_{n}\left(z\right)}.
\]
A necessary and sufficient condition for the point at infinity to
be at most a regular singularity of this operator is that the function
$r_{j}\left(z\right)$ is analytic at $\infty$ and has a zero there
of order at least $j$ (\cite[Theorem 9.8b]{henrici1977applied}).
That is, 
\[
\deg p_{n}-\deg p_{n-j}\geq j.
\]
But this is equivalent to
\begin{eqnarray*}
\deg p_{n}-n & \geqslant & \deg p_{n-j}-\left(n-j\right)\\
\alpha_{n} & \geqslant & \alpha_{n-j}.
\end{eqnarray*}
To prove the second statement, substitute the formal Frobenius series
at infinity
\[
g\left(z\right)=\sum_{k=0}^{\infty}\frac{b_{k}}{z^{s+k}}
\]
into $\Op g=0$. By complete analogy with the calculation in \prettyref{thm:cauchy-int-ode}
we get the recurrence relation
\[
\sum_{j=0}^{n}\sum_{i=0}^{d_{j}}\left(-1\right)^{j}a_{i,j}\ff{t+s+i-1}jb_{t+i-j}=0,\qquad t=0,1,\dots.
\]
For $t=0$ we find the highest order coefficient in this recurrence
to be equal to ($i-j=\alpha=\alpha_{n}$)
\[
\sum_{j=0}^{n}\left(-1\right)^{j}\ff{s+\alpha+j-1}ja_{j+\alpha,j}=q_{\alpha}\left(s-1\right).
\]
The proof is finished.
\end{svmultproof2}
Together with \prettyref{thm:mvi-bound-for-all}, this immediately
implies the following bound.
\begin{corollary}
Let $\Op\in\prp$, and let $\lambda\left(\Op\right)$ denote its largest
positive integer characteristic exponent at the point $\infty$. Then
$\Lambda\left(\Op\right)=\lambda\left(\Op\right)-1$, and consequently
\[
\mvi\left(\Op,\np\right)\leqslant\max\left\{ \left(\np+2\right)\ord,\lambda\left(\Op\right)\right\} +\dg_{n}-n-1.
\]

\end{corollary}
Now let us return to \prettyref{ex:legendre}. The following fact
is well-known (e.g. \cite[Section 9.10]{henrici1977applied}.
\begin{proposition}
The Legendre differential operator $\Op_{m}$ is of Fuchsian type
with singularities $-1,1,\infty$. The characteristic exponents at
$\infty$ are $m+1$ and $-m$.
\end{proposition}
\prettyref{thm:mvi-general-case} is clearly not applicable. Using
the formula \eqref{eq:ek-def}, it is easy to see that 
\[
P_{\Op}\left(\fun,x^{k}\right)\left(1\right)=P_{\Op}\left(\fun,x^{k}\right)\left(-1\right)=0
\]
for any $\fun\in\PD$, and therefore the sequence $\left\{ \e_{k}\right\} $
in this case is identically zero. Consequently, we conclude that
\[
\mvi\left(\Op_{m},0\right)=m,
\]
as expected.

\bibliographystyle{plain}
\bibliography{../../../bibliography/all-bib}

\begin{thebibliography}{10}

\bibitem{ang2002mta}
D.D. Ang.
\newblock {\em {Moment Theory and Some Inverse Problems in Potential Theory and
  Heat Conduction}}.
\newblock Springer, 2002.

\bibitem{batFullFourier}
D.~Batenkov.
\newblock {Complete Algebraic Reconstruction of Piecewise-Smooth Functions from
  Fourier Data}.
\newblock {\em arXiv preprint arXiv:1211.0680}.

\bibitem{bat2008}
D.~Batenkov.
\newblock {Moment inversion problem for piecewise D-finite functions}.
\newblock {\em Inverse Problems}, 25(10):105001, October 2009.

\bibitem{BatGolYom2011}
D.~Batenkov, V.~Golubyatnikov, and Y.~Yomdin.
\newblock {Reconstruction of Planar Domains from Partial Integral
  Measurements}.
\newblock In {\em Proc. Complex Analysis \& Dynamical Systems V}, 2011.

\bibitem{batyomTaylrPoincBarc}
D.~Batenkov and Y.~Yomdin.
\newblock {Taylor Domination, Tur\'an lemma, and Poincar\'e-Perron Sequences}.
\newblock {\em Submitted to this volume.}

\bibitem{batyomAlgFourier}
D.~{Batenkov} and Y.~{Yomdin}.
\newblock {Algebraic Fourier reconstruction of piecewise smooth functions}.
\newblock {\em Mathematics of Computation}, 81:277--318, 2012.

\bibitem{briskin2010center}
M.~Briskin, N.~Roytvarf, and Y.~Yomdin.
\newblock {Center conditions at infinity for Abel differential equations}.
\newblock {\em {Annals of Mathematics}}, 172(1):437--483, 2010.

\bibitem{elaydi2005ide}
S.~Elaydi.
\newblock {\em {An Introduction to Difference Equations}}.
\newblock Springer, 2005.

\bibitem{gustafsson2000rpd}
B.~Gustafsson, C.~He, P.~Milanfar, and M.~Putinar.
\newblock {Reconstructing planar domains from their moments}.
\newblock {\em Inverse Problems}, 16(4):1053--1070, 2000.

\bibitem{henrici1977applied}
P.~Henrici.
\newblock {\em Applied and Computational Complex Analysis: Vol.: 2.: Special
  Functions: Integral Transforms: Asymptotics: Continued Fractions}.
\newblock John Wiley \& Sons, 1977.

\bibitem{ince1956ode}
E.L. Ince.
\newblock {\em {Ordinary Differential Equations}}.
\newblock Courier Dover Publications, 1956.

\bibitem{kisunko1}
V.~Kisunko.
\newblock {Cauchy Type Integrals and a D-moment Problem}.
\newblock {\em Mathematical Reports of the Academy of Science of the Royal
  Society of Canada}, 29(4):115--122, 2008.

\bibitem{myerPort1995}
G.~Myerson and A.~J. van~der Poorten.
\newblock Some problems concerning recurrence sequences.
\newblock {\em The American Mathematical Monthly}, 102(8):pp. 698--705, 1995.

\bibitem{pakovich:cti}
F.~Pakovich, N.~Roytvarf, and Y.~Yomdin.
\newblock Cauchy-type integrals of algebraic functions.
\newblock {\em Israel Journal of Mathematics}, 144(2):221--291, 2004.

\bibitem{vetterli2002ssf}
M.~Vetterli, P.~Marziliano, and T.~Blu.
\newblock {Sampling signals with finite rate of innovation}.
\newblock {\em IEEE Transactions on Signal Processing}, 50(6):1417--1428, 2002.

\end{thebibliography}

\end{document}